\documentclass[11pt,a4paper]{article}

\usepackage{amsmath}
\usepackage{amsfonts}

\usepackage[latin1]{inputenc}
\usepackage{fontenc}
\usepackage{euscript}
\usepackage{array}
\usepackage{vmargin}
\usepackage{amssymb}
\usepackage{latexsym}
\usepackage{enumerate}
\usepackage{multirow}
\usepackage{graphicx}

\newcommand{\pg}{\hspace{0.6cm}}

\def\mathbi#1{\textbf{\em #1}}

\def\N{{\rm I\kern-.20em N}}
\def\R{{\rm I\kern-.20em R}}
\def\indi{{1\kern-.20em\rm I}}

\newtheorem{definition}{Definition}[section]
\newtheorem{proposition}{Proposition}[section]
\newtheorem{theorem}{Theorem}
\newtheorem{remark}{Remark}

\newcommand {\edem}{\hfill $\square$ \end {proof}}
\newtheorem{ex}{Example}[section]
\newcommand{\bdem} {\begin{proof}}

\begin{document}
\bibliographystyle{plain}
\hyphenation{o-pe-ra-tors}

\title{Generalized madogram and pairwise dependence of maxima over two regions of a random field}
\author{Fonseca, C.\footnote{\noindent Instituto Politécnico da Guarda, Portugal, E-mail: {\tt
cfonseca@ipg.pt}\hspace*{2mm}},  Pereira, L. and
Ferreira, H. and Martins, A.P. \footnote{\noindent Departamento de Matemática,
Universidade da Beira Interior, Portugal, E-mail: {\tt
lpereira@ubi.pt}, {\tt helenaf@ubi.pt}, {\tt
amartins@ubi.pt}\hspace*{2mm}}}


\date{}
\maketitle

\noindent {\bf Abstract:} Spatial environmental processes often exhibit dependence in their large values. In order to model such processes their dependence properties must be characterized and quantified. In this paper we introduce a measure that evaluates the dependence among extreme observations located in two separated regions of locations of $\mathbb{R}^2$. We compute the range of this new dependence measure, which extends the existing $\lambda$-madogram concept, and compare it with extremal coefficients, finding generalizations of the known relations in pairwise approach. Estimators for this measure are introduced and asymptotic normality and strong consistency are shown. An application to the annual maxima precipitation in Portuguese regions is presented.





\pagenumbering{arabic}

\section{Introduction}

\pg Natural models for spatial extremes, as observed in environmental, atmospheric and geological sciences, are max-stable processes. These processes arise from an infinite-dimensional generalization of extreme value theory and date back to de Haan (1984) \cite{haan1}, Vatan (1985) \cite{vat} and de Haan and Pickands (1986) \cite{haan2}, who obtained, among other results, a spectral representation of such processes.

Max-stable processes can be, for example, good approximations for annual maxima of daily spatial rainfall (Smith \cite{smi}, Coles \cite{col}, Schlather \cite{sch1}, among others) and therefore have been widely applied to real data.

Briefly, a max-stable process ${\bf Z}=\left\{Z_{\mathbi{x}}\right\}_{{\mathbi{x}}\in \mathbb R^d}$ is the limit process of maxima of i.i.d. random fields $Y^{(i)}_{\mathbi{x}},$ ${\mathbi{x}}\in \R^d,$ $i=1,\ldots,n.$ Namely, for suitable $a_n({\mathbi{x}})>0$ and $b_n({\mathbi{x}})\in \R,$ $$Z_{{\mathbi{x}}}=\lim_{n\to \infty}\frac{\bigvee_{i=1}^nY^{(i)}_{\mathbi{x}}-b_n({\mathbi{x}})}{a_n({\mathbi{x}})},\qquad {\mathbi{x}}\in \R^d,$$ provided the limit exists.

We shall consider $d=2,$ that is ${\bf Z}=\left\{Z_{\mathbi{x}}\right\}_{{\mathbi{x}}\in \mathbb R^2}.$ The distribution of $(Z_{\mathbi{x}_1},\ldots,Z_{\mathbi{x}_k})$ is a multivariate extreme value (MEV) distribution $G$ and we can assume, without loss of generality, that the margins of ${\bf Z}$ have a unit Fréchet distribution, $F(x)=\exp(-x^{-1}),\ x>0$ (Resnick \cite{res}). The distribution $G$ can then be defined by $$G(z_1,\ldots,z_k)=\exp(-V_{\bf{z}}(z_1,\ldots,z_k)),\quad z_i\in \R^+,$$ where $V_{\bf{z}}$ denotes the dependence function of the MEV distribution $G$, which is homogeneous of order $-1,$ i.e., $V_{\bf{z}}(\alpha z_1,\ldots,\alpha z_k)=\alpha^{-1}V_{\bf{z}}(z_1,\ldots,z_k),$ $ z_i\in \R^+,$ $\alpha>0.$

The dependence function captures the multivariate dependence structure and the scalar $V_{\bf{z}}(1,\ldots,1)$  defines the extremal coefficient considered in Schlather an Tawn \cite{sch} which measures the extremal dependence between the variables indexed in the region ${\bf{z}}$. This coefficient varies between 1 and $k$ depending on the degree of dependence among the $k$ variables. These measures of dependence have gained great importance since quantifying dependence between extreme events occurring at several locations of a random field is a fundamental issue in applied spatial extreme value analysis.

Cooley {\it{et al.}} (2006) showed that the bivariate extremal coefficient can be directly estimated from the madogram (that represents a first order variogram), obtaining in this way a connection between extreme value theory and the field of geostatistics.  A measure of the full pairwise extremal dependence function is given by the $\lambda$-madogram defined in Naveau \textit{et al}. \cite{na} as
\begin{equation}
\nu^{\lambda}(\mathbi{x}_1,\mathbi{x}_2)=\frac{1}{2}E\left|F^{\lambda}(Z_{\mathbi{x}_1})-F^{1-\lambda}(Z_{\mathbi{x}_2})\right|,\quad \lambda\in (0,1),\label{lambda}
\end{equation}
where $F$ denotes the marginal distribution of ${\bf Z}.$

Although the $\lambda$-madogram fully characterizes the pairwise extremal dependence it does not enable the analysis of dependence between maxima over two disjoint regions of locations.  The importance of characterizing dependence between extremes occurring at two disjoint regions has been recognized by hydrologists who have grouped data into regions based on geographical or catchment characteristics. It is clear that, for example, the weather in mountain regions usually affects the weather in the surrounding regions and therefore a measure able to capture such regional dependence is essential.

In this paper we propose a measure that enables the analysis of dependence between maxima over two disjoint regions of locations $\textbf{x}=\{\mathbi{x}_1,\ldots,\mathbi{x}_k\}$ and $\textbf{y}=\{\mathbi{y}_1,\ldots,\mathbi{y}_s\}$ and therefore generalizes the $\lambda$-madogram. This measure, here called generalized madogram, is introduced in Section 2 and some of its main properties are presented, namely its relation with the dependence function of the MEV G. In Section 3 we present estimators for the generalized madogram and derive the respective properties of strong consistency and asymptotic normality. The performance of the proposed estimators is analyzed in Section 4 with a max-stable M4 random field. Finally, Section 5 illustrates our approach through an application to precipitation data from Portugal.

\section{Generalized madogram and dependence of spatial extreme events}

\pg The measure that enables the analysis of dependence between maxima over two disjoint regions of locations  ${\textbf{x}}=\{\mathbi{x}_1,\ldots,\mathbi{x}_k\}$ and $\textbf{y}=\{\mathbi{y}_1,\ldots,\mathbi{y}_s\}$ and generalizes the $\lambda$-madogram introduced in Naveau \textit{et al}. \cite{na} is defined as follows.

 \begin{definition}
 Let $Z = \{Z_{\textbf{x}}\}_{{\textbf{x}}\in\R^2}$ be a max-stable random field with unit Fréchet margins and ${\mathbi{x}}=\{\textbf{x}_1, \ldots, \textbf{x}_k\}$ and ${\mathbi{y}} =\{\textbf{y}_1,\ldots,\textbf{y}_s\}$ two disjoint regions of $\R^2.$ The generalized madogram is defined as
\begin{equation}
\nu^{\alpha,\beta}({\bf{x}},{\bf{y}})=\frac{1}{2}E\left|F^{\alpha}(M({\bf{x}}))-F^{\beta}(M({\bf{y}}))\right|,\quad \alpha> 0,\  \beta> 0 ,\label{mad_gen}
\end{equation}
where  $M({\mathbi{x}})=\bigvee_{i=1}^k Z_{{\it{\textbf{x}}}_i}$ and $M({\mathbi{y}})=\bigvee_{j=1}^s Z_{\textbf{y}_j}$.
\end{definition}\vspace{0.3cm}

\begin{remark}
{\rm{The following equalities hold for the generalized madogram
\begin{eqnarray*}
\nu^{\alpha,\beta}(\textbf{x},\textbf{y})&=&\frac{1}{2}E\left|\bigvee_{i=1}^k F^{\alpha}(Z_{\mathbi{x}_i})-\bigvee_{j=1}^s F^{\beta}(Z_{\mathbi{y}_j})\right|\\[0.3cm]
&=&\frac{1}{2}E\left|\bigvee_{i=1}^k F \left( \frac{Z_{\mathbi{x}_i}}{\alpha}\right)-\bigvee_{j=1}^s F\left(\frac{Z_{\mathbi{y}_j}}{\beta}\right)\right|.
\end{eqnarray*}\vspace{0.3cm}

\noindent This representation of the generalized madogram, $\nu^{\alpha,\beta}(\textbf{x},\textbf{y}),$ will motivate the definition of natural estimators, in Section 3, for this coefficient.}}
\end{remark}

\begin{remark}
{\rm{When we take $\beta=1-\alpha$, $\alpha\in (0,1)$, and $k=s=1$ in (\ref{mad_gen}), we obtain (\ref{lambda}).}}
\end{remark}

The following proposition states that $\nu^{\alpha,\beta}(\textbf{x},\textbf{y})$ provides dependence information between the regions $\mathbf{x}$ and $\mathbf{y}$ through the dependence function of the MEV distribution $G.$ This result generalizes Proposition 1. in Naveaux \emph{et al.} \cite{na}.

\begin{proposition}
For any max-stable random field with unit Fréchet margins and for each pair of disjoint regions of locations ${\bf{x}}=\{{\textbf{x}}_1,\ldots,{\textbf{x}}_k\}$ and ${\bf{y}}=\{{\textbf{y}}_1,\ldots,{\textbf{y}}_s\}$ in $\R^2$, we have
$$
\nu^{\alpha,\beta}({\bf x , y} )=\frac{V_{{\bf x ,y}}(\alpha,\ldots,\alpha,\beta,\ldots,\beta)}{1+V_{\mathbf{x},\mathbf{y}}
(\alpha,\ldots,\alpha,\beta,\ldots,\beta)}-c(\alpha,\beta)$$ with
$$c(\alpha,\beta)=\frac{1}{2}\left(\frac{V_{\mathbf{x}}(1,\ldots,1)}{\alpha+V_{\mathbf{x}}(1,\ldots,1)}+\frac{V_{\mathbf{y}}(1,\ldots,1)}
{\beta+V_{\mathbf{y}}(1,\ldots,1)}\right),
$$
where $$V_{\mathbf{x},\mathbf{y}}(z_1,\ldots,z_k,z_{k+1},\ldots,z_{k+s})=-\ln G_{\mathbf{x},\mathbf{y}}(z_1,\ldots,z_k,z_{k+1},\ldots,z_{k+s})$$ and $$G_{\mathbf{x},\mathbf{y}}(z_1,\ldots,z_{k+s})=P\left(\left\{\bigcap_{i=1}^k\left\{Z_{{\textbf{y}}_i}\leq z_i\right\}\right\}\bigcap\left\{\bigcap_{i=1}^s\left\{Z_{{\textbf{y}}_i}\leq z_{k+i}\right\}\right\}\right), \quad z_i\in \R^+.$$
\end{proposition}

\bdem To obtain the result, we start by transforming the definition of $\nu^{\alpha,\beta}(\textbf{x},\textbf{y})$ through the relation \linebreak $|a-b|=2(a \vee b)-(a+b)$, and then take into account that
$$
E\left(F^\alpha\left(M(\bf x)\right)\vee F^\beta\left(M(\bf y)\right)\right)=E\left(F\left(\frac{M(\bf x)}{\alpha}\vee \frac{M(\bf y)}{\beta}\right)\right).
$$
It holds 
\begin{eqnarray*}
P\left(\frac{M(\bf x)}{\alpha}\vee \frac{M(\bf y)}{\beta}\leq u\right)&=&P\left(M({\bf x})\leq \alpha u, M({\bf y})\leq \beta u\right) \\[0.3cm]
&=&G_{\bf x, \bf y}(\alpha u,\ldots,\alpha u,\beta u,\ldots,\beta u)\\[0.3cm]
&=&\exp\left\{-u^{-1}V_{\bf x, \bf y}(\alpha,\ldots,\alpha,\beta,\ldots,\beta)\right\}, \ \ u>0.
\end{eqnarray*}
Hence,
\begin{eqnarray*}
\lefteqn{E\left(F\left(\frac{M(\bf x)}{\alpha}\vee \frac{M(\bf y)}{\beta}\right)\right)}\\[0.3cm]
&=&\int_0^{+\infty}F(u)\exp\left(-V_{\bf x,\bf y}(\alpha u,\ldots,\alpha u,\beta u,\ldots,\beta u)\right)\frac{d}{du}\left(-V_{\bf x,\bf y}(\alpha u,\ldots,\alpha u,\beta u,\ldots,\beta u)\right)\\[0.3cm]
&=&\int_0^{+\infty}\exp\left(-u^{-1}-u^{-1}V_{\bf x,\bf y}(\alpha ,\ldots,\alpha ,\beta ,\ldots,\beta )\right)u^{-2}V_{\bf x,\bf y}(\alpha ,\ldots,\alpha ,\beta ,\ldots,\beta )du\\[0.3cm]
&=&\frac{V_{\bf x,\bf y}(\alpha,\ldots,\alpha,\beta,\ldots,\beta)}{1+V_{\bf x,\bf y}(\alpha,\ldots,\alpha,\beta,\ldots,\beta)}
\end{eqnarray*}

Using similar arguments, we obtain
$E\left(F^\alpha(M(\bf x))\right)=\frac{V_{\bf x}(\alpha,\ldots,\alpha)}{1+V_{\bf x}(\alpha,\ldots,\alpha)}=\frac{V_{\bf x}(1,\ldots,1)}{\alpha+V_{\bf x}(1,\ldots,1)}
$ and\linebreak  $E\left(F^\beta(M(\bf y))\right)=\frac{V_{\bf y}(\beta,\ldots,\beta)}{1+V_{\bf y}(\beta,\ldots,\beta)}=\frac{V_{\bf y}(1,\ldots,1)}{\beta+V_{\bf y}(1,\ldots,1)}$.
\edem\vspace{0.3cm}

\begin{remark}{\rm{
 For each $\alpha,\ \beta> 0$ the coefficient $c(\alpha, \beta)$ considers the dependence intra each of the regions ${\mathbf{x}}$ and ${\mathbf{y}}$ through the extremal coefficients of vectors with margins $Z_{\mathbi{x}_1},\ldots,Z_{\mathbi{x}_k}$ and $Z_{\mathbi{y}_1},\ldots,Z_{\mathbi{y}_s}$. When we consider $c(\alpha, \beta)$ constant, the dependence between ${\mathbf{x}}$ and ${\mathbf{y}}$ is stronger for lower values of $\nu^{\alpha,\beta}(\textbf{x},\textbf{y})$, corresponding to lower values of $V_{{\bf x ,y}}(\alpha,\ldots,\alpha,\beta,\ldots,\beta)$.}}
\end{remark}\vspace{0.3cm}


In the following proposition we establish some properties of the generalized madogram.

\begin{proposition}
Let ${\mathbf{x}}=\{{\textbf{x}}_1,\ldots,{\textbf{x}}_k\}$ and ${\mathbf{y}}=\{{\textbf{y}}_1,\ldots,{\textbf{y}}_s\}$ be disjoint regions of $\R^2$. We have, for each $\alpha,\beta\in \R^+$,
\begin{enumerate}
\item[\bf 1.] $0 \leq \nu^{\alpha,\beta}(\mathbf{x},\mathbf{y})\leq \frac{1}{2};$
\item[\bf 2.] $\nu^{\alpha,\alpha}(\mathbf{x},\mathbf{y})=\frac{\epsilon_{\bf x \cup \bf y}}{\alpha+\epsilon_{\bf x\cup \bf y}}-
\frac{1}{2}\left(\frac{\epsilon_{\mathbf{x}}}{\alpha+\epsilon_{\mathbf{x}}}+\frac{\epsilon_{\mathbf{y}}}{\alpha+\epsilon_{\mathbf{y}}}\right),$ \ \ where $\epsilon_{\bf x \cup \bf y}=V_{\mathbf{x},\mathbf{y}}(1,\ldots,1)$.
\end{enumerate}
\end{proposition}\vspace{0.3cm}

\bdem  The first statement results from the definition of the generalized madogram and the second follows from the definition of the extremal coefficient $\epsilon$.
\edem\vspace{0.3cm}

\begin{remark}
{\rm{The function $\nu^{\alpha,\alpha}(\mathbf{x},\mathbf{y})$ can also be related with the dependence coefficients considered in Ferreira \cite{fer} as follows:
$$\nu^{\alpha,\alpha}(\mathbf{x},\mathbf{y})=\displaystyle{\frac{\epsilon_{\mathbf{y}}\epsilon_1(\mathbf{x},\mathbf{y})}
{\alpha+\epsilon_{\mathbf{y}}\epsilon_1(\mathbf{x},\mathbf{y})}}-c(\alpha,\alpha)=\displaystyle{\frac{(\epsilon_{\mathbf{y}}+\epsilon_{\mathbf{x}})
\epsilon_2(\mathbf{x},\mathbf{y})}{\alpha+(\epsilon_{\mathbf{y}}+\epsilon_{\mathbf{x}})\epsilon_2(\mathbf{x},\mathbf{y})}}-
c(\alpha,\alpha),\;\;\alpha>0, $$
where
$\epsilon_1({\mathbf{x}},{\mathbf{y}})=\displaystyle{\frac{\epsilon_{\bf x \cup \bf y}}
{\epsilon_{\mathbf{y}}}}$ and
$\epsilon_2({\mathbf{x}},{\mathbf{y}})=\displaystyle{\frac{\epsilon_{\bf x \cup \bf y}}
{\epsilon_{\mathbf{x}}+\epsilon_{\mathbf{y}}}}.$ These coefficients evaluate the strength of dependence between the events $\left\{M({\mathbf{x}})\leq u\right\}$ and $\left\{M({\mathbf{y}})\leq u\right\}.$
}}
\end{remark}\vspace{0.3cm}

\begin{remark}
{\rm{If the variables $M(\mathbf{x})$ and $M(\mathbf{y})$ are independent then
$$
\nu^{\alpha,\alpha}(\mathbf{x},\mathbf{y})=\frac{\epsilon_{\mathbf{x}}+\epsilon_{\mathbf{y}}}{\alpha+\epsilon_{\mathbf{x}}+\epsilon_{\mathbf{y}}}-\frac{1}{2}\left(\frac{\epsilon_{\mathbf{x}}}{\alpha+\epsilon_{\mathbf{x}}}+\frac{\epsilon_{\mathbf{y}}}{\alpha+\epsilon_{\mathbf{y}}}\right),
$$ whereas if  they are totally dependent
$$\nu^{\alpha,\alpha}(\mathbf{x},\mathbf{y})=\frac{\epsilon_{\mathbf{x}}\vee\epsilon_{\mathbf{y}}}{\alpha+\epsilon_{\mathbf{x}}\vee\epsilon_{\mathbf{y}}}-\frac{1}{2}\left(\frac{\epsilon_{\mathbf{x}}}{\alpha+\epsilon_{\mathbf{x}}}+\frac{\epsilon_{\mathbf{y}}}{\alpha+\epsilon_{\mathbf{y}}}\right).
$$}}
\end{remark}


\section{Estimating the generalized madogram}

\pg Proposition 2.1 relates the generalized madogram with well known dependence measures. Immediate estimators for the generalized madogram can then be obtained through the estimators of those measures, which have already been studied in the literature. For a survey we refer to Krajina \cite{kraj}, Beirlant \cite{beir} and Schlather and Tawn \cite{sch}.

In this section we present a natural non-parametric estimator for the generalized madogram based on sample means.\vspace{0.3cm}

Let $(Z_{{\mathbi{x}}_1}^{(t)},\ldots,Z_{{\mathbi{x}}_k}^{(t)})$ and $(Z_{\mathbi{y}_1}^{(t)},\ldots,Z_{\mathbi{y}_s}^{(t)}),$ $t=1,\ldots,T,$ be independent replications of $(Z_{\mathbi{x}_1},\ldots,Z_{\mathbi{x}_k})$ and $(Z_{\mathbi{y}_1},\ldots,Z_{\mathbi{y}_s}),$ respectively. Hence $\{M_t({\mathbf{x}})=\bigvee_{i=1}^k Z_{{\mathbi{x}}_i}^{(t)},\ t=1,\ldots, T\}$ and $\{M_t({\mathbf{y}})=\bigvee_{i=1}^s Z_{{\mathbi{y}}_i}^{(t)},$ $t=1,\ldots, T\}$ are random samples of $M(\mathbf{x})$ and $M(\mathbf{y}),$ respectively.

In the case of known maginal distribution $F$, which becomes unit Fréchet by transformation $-\frac{1}{\log F(Z_{\mathbi{x}_i})}$ for $\mathbi{x}_i\in \mathbb R^2$, the estimator for the generalized madogram is given by $$\widehat{\nu}^{\alpha,\beta}(\mathbf{x},\mathbf{y})=\frac{1}{2T}\sum_{i=1}^T|F^{\alpha}(M_i({\mathbf{x}}))-F^{\beta}(M_i({\mathbf{y}}))|,\quad \alpha > 0,\ \beta > 0.$$
This estimator is unbiased and converges in distribution to a Gaussian distribution, as stated in the following proposition.
\begin{proposition} (Asymptotic normality and strong consistency under known marginal distribution F)

 We have
 $$
 \frac{\sqrt{T}(\widehat{\nu}^{\alpha,\beta}(\mathbf{x},\mathbf{y})-\nu^{\alpha,\beta}(\mathbf{x},\mathbf{y}))}{\sigma}\rightarrow N(0,1),
 $$
where $\sigma^2=\frac{1}{2}\gamma_F^{\alpha,\beta}(\mathbf{x},\mathbf{y})-(\nu^{\alpha,\beta}(\mathbf{x},\mathbf{y}))^2$ and $\gamma_F^{\alpha,\beta}(\mathbf{x},\mathbf{y})=\frac{1}{2}E\left[(F^\alpha(M(\mathbf x))-F^\beta(M(\mathbf y)))^2\right]$. Moreover,  $\widehat{\nu}^{\alpha,\beta}(\mathbf{x},\mathbf{y})$ converges almost surely to ${\nu}^{\alpha,\beta}(\mathbf{x},\mathbf{y})$.
\end{proposition}
\bdem
Let $Y_1,\ldots,Y_T$ be independent copies of $Y=\frac{1}{2}\left|F^\alpha(M(\mathbf x))-F^\beta(M(\mathbf y))\right|$. We have that
$$
\frac{\sqrt{T}(\bar{Y}-\mu_Y)}{\sigma_Y}\rightarrow N(0,1),
$$
where $\mu_Y=\frac{1}{2}E\left|F^\alpha(M(\mathbf x))-F^\beta(M(\mathbf y))\right|=\nu^{\alpha,\beta}(\mathbf{x},\mathbf{y})$ and $\sigma^2_Y=\frac{1}{2}\gamma_F^{\alpha,\beta}(\mathbf{x},\mathbf{y})-(\nu^{\alpha,\beta}(\mathbf{x},\mathbf{y}))^2$.

The strong consistency of $\widehat{\nu}^{\alpha,\beta}(\mathbf{x},\mathbf{y})$ follows since the sample mean converges almost surely to the mean value.\edem

When the distribution of each $Z_{\mathbi{x}_i}$, $F_{\mathbi{x}_i}$, is unknown we take the empirical Fréchet normalization of the variables, i.e., $\widehat{U}_{\mathbi{x}_i}^{(t)}=-\frac{1}{\log(\widehat{F}_{\mathbi{x}_i}(Z^{(t)}_{\mathbi{x}_i}))},$ where $\widehat{F}_{\mathbi{x}_i}(Z^{(t)}_{\mathbi{x}_i})$ is the empirical distribution function. We then find the following modification of the above estimator
\begin{eqnarray*}
\widehat{\widehat{\nu}}^{\alpha,\beta}(\mathbf{x},\mathbf{y})& = & \frac{1}{2T}\sum_{t=1}^T \left|F^{\alpha}\left(\bigvee_{i=1}^k \widehat{U}^{(t)}_{{\mathbi{x}}_i}\right)-
F^{\beta}\left(\bigvee_{j=1}^s \widehat{U}^{(t)}_{{\mathbi{y}}_j}\right)\right|\\[0.3cm]
& = & \frac{1}{2T}\sum_{t=1}^T\left|\bigvee_{i=1}^k F\left( \frac{\widehat{U}^{(t)}_{{\mathbi{x}}_i}}{\alpha}\right)-
\bigvee_{j=1}^s F\left(\frac{\widehat{U}^{(t)}_{{\mathbi{y}}_j}}{\beta}\right)\right| \\[0.3cm]
& = & \frac{1}{2T}\sum_{t=1}^T\left|\bigvee_{i=1}^k \widehat{F}^{\alpha}_{{\mathbi{x}}_i}\left( Z^{(t)}_{{\mathbi{x}}_i}\right)-
\bigvee_{j=1}^s\widehat{F}^{\beta}_{{\mathbi{y}}_j}\left(Z^{(t)}_{{\mathbi{y}}_j}\right)\right| ,\quad \alpha> 0,\ \beta> 0,\label{estim}
\end{eqnarray*}
where $$\widehat{F}_{{\mathbi{x}}_i}(u)=\frac {1}{T}\sum_{t=1}^T \indi_{\{Z_{{\mathbi{x}}_i}^{(t)}\leq u\}}.$$\vspace{0.3cm}




This estimator, $\widehat{\widehat{\nu}}^{\alpha,\beta}(\mathbf{x},\mathbf{y})$ is strongly consistent since it holds
\begin{eqnarray*}
\lefteqn{\left|\frac{1}{2}\frac{1}{T}\sum_{t=1}^T\left|\bigvee_{i=1}^k \widehat{F}^{\alpha}_{{\mathbi{x}}_i}\left( Z^{(t)}_{{\mathbi{x}}_i}\right)-
\bigvee_{j=1}^s\widehat{F}^{\beta}_{{\mathbi{y}}_j}(Z^{(t)}_{{\mathbi{y}}_j})\right|-\frac{1}{2}E\left|\bigvee_{i=1}^k F^\alpha(Z^{(t)}_{{\mathbi{y}}_j})-\bigvee_{j=1}^s F^{\beta}(Z^{(t)}_{{\mathbi{y}}_j})\right|\right|}\\[0.3cm]
&\leq& \left|\frac{1}{2}\frac{1}{T}\sum_{t=1}^T\left|\bigvee_{i=1}^k \widehat{F}^{\alpha}_{{\mathbi{x}}_i}\left( Z^{(t)}_{{\mathbi{x}}_i}\right)-
\bigvee_{j=1}^s\widehat{F}^{\beta}_{{\mathbi{y}}_j}(Z^{(t)}_{{\mathbi{y}}_j})\right|-\frac{1}{2}\frac{1}{T}\sum_{t=1}^T\left|\bigvee_{i=1}^k F^\alpha(Z^{(t)}_{{\mathbi{x}}_i})-\bigvee_{j=1}^s F^{\beta}(Z^{(t)}_{{\mathbi{y}}_j})\right|\right|\\[0.3cm]
&&+\left|\frac{1}{2}\frac{1}{T}\sum_{t=1}^TE\left|\bigvee_{i=1}^k F^\alpha(Z^{(t)}_{{\mathbi{x}}_i})-\bigvee_{j=1}^s F^{\beta}(Z^{(t)}_{{\mathbi{y}}_j})\right|-\frac{1}{2}E\left|\bigvee_{i=1}^k F^\alpha(Z^{(t)}_{{\mathbi{x}}_i})-\bigvee_{j=1}^s F^{\beta}(Z^{(t)}_{{\mathbi{y}}_j})\right|\right|.
\end{eqnarray*}
The second term converges almost surely to zero by the strong Law of Large Numbers.

In what concerns the first term we have
\begin{eqnarray*}
\lefteqn{\left|\frac{1}{2}\frac{1}{T}\sum_{t=1}^T\left|\bigvee_{i=1}^k \widehat{F}^{\alpha}_{{\mathbi{x}}_i}\left( Z^{(t)}_{{\mathbi{x}}_i}\right)-
\bigvee_{j=1}^s\widehat{F}^{\beta}_{{\mathbi{y}}_j}(Z^{(t)}_{{\mathbi{y}}_j})\right|-\frac{1}{2}\frac{1}{T}\sum_{t=1}^T\left|\bigvee_{i=1}^k F^\alpha(Z^{(t)}_{{\mathbi{x}}_i})-\bigvee_{j=1}^s F^{\beta}(Z^{(t)}_{{\mathbi{y}}_j})\right|\right|}\\[0.3cm]
&\leq&\frac{1}{2}\frac{1}{T}\sum_{t=1}^T \left| \left|\bigvee_{i=1}^k \widehat{F}^{\alpha}_{{\mathbi{x}}_i}\left( Z^{(t)}_{{\mathbi{x}}_i}\right)-
\bigvee_{j=1}^s\widehat{F}^{\beta}_{{\mathbi{y}}_j}(Z^{(t)}_{{\mathbi{y}}_j})\right|-\left|\bigvee_{i=1}^k F^\alpha(Z^{(t)}_{{\mathbi{x}}_i})-\bigvee_{j=1}^s F^{\beta}(Z^{(t)}_{{\mathbi{y}}_j})\right|\right|\\[0.3cm]
&\leq&\frac{1}{2}\frac{1}{T}\sum_{t=1}^T \left| \bigvee_{i=1}^k \widehat{F}^{\alpha}_{{\mathbi{x}}_i}\left( Z^{(t)}_{{\mathbi{x}}_i}\right)-
\bigvee_{j=1}^s\widehat{F}^{\beta}_{{\mathbi{y}}_j}(Z^{(t)}_{{\mathbi{y}}_j})-\bigvee_{i=1}^k F^\alpha(Z^{(t)}_{{\mathbi{x}}_i})+\bigvee_{j=1}^s F^{\beta}(Z^{(t)}_{{\mathbi{y}}_j})\right|\\[0.3cm]
&\leq&\frac{1}{2}\frac{1}{T}\sum_{t=1}^T\left[\bigvee_{i=1}^k \left|\widehat{F}^{\alpha}_{{\mathbi{x}}_i}\left( Z^{(t)}_{{\mathbi{x}}_i}\right)-
F^\alpha(Z^{(t)}_{{\mathbi{x}}_i})\right|+\bigvee_{j=1}^s\left|\widehat{F}^{\beta}_{{\mathbi{y}}_j}(Z^{(t)}_{{\mathbi{y}}_j})- F^{\beta}(Z^{(t)}_{{\mathbi{y}}_j})\right|\right]\\[0.3cm]
&\leq&\frac{1}{2}\frac{1}{T}\sum_{t=1}^T\left[\sum_{i=1}^k \left|\widehat{F}^{\alpha}_{{\mathbi{x}}_i}\left( Z^{(t)}_{{\mathbi{x}}_i}\right)-
F^\alpha(Z^{(t)}_{{\mathbi{x}}_i})\right|+\sum_{j=1}^s\left|\widehat{F}^{\beta}_{{\mathbi{y}}_j}(Z^{(t)}_{{\mathbi{y}}_j})- F^{\beta}(Z^{(t)}_{{\mathbi{y}}_j})\right|\right],
\end{eqnarray*}
which converges almost surely to zero. This follows from the strong consistency of the empirical distribution function and the fact that $\frac{1}{T}\sum_{i=1}^T|X_i|$ converges almost surely to $E|X_1|$, where $X_i$, $i=1,\ldots,T$, are i.i.d. random variables with $E|X_1|<\infty$.\\
\smallskip \\
The asymptotic normality of the estimator is obtained by taking $J(x_1,\ldots,x_{k+s})=\frac{1}{2}\left|\bigvee_{i=1}^kx_i^{\alpha}-\bigvee_{i=k+1}^{k+s}x_i^{\beta}\right|$ in the following theorem stated in Fermanian \textit{et al.} \cite{ferman}.

\begin{theorem}
Let $(Z_{\textbf{x}_1},\ldots,Z_{\textbf{x}_k},Z_{\textbf{y}_1},\ldots,Z_{\textbf{y}_s})$ be a random vector with d.f. $H$ and continuous marginal d.f.'s $F_{\textbf{x}_1},\ldots,F_{\textbf{x}_k}, F_{\textbf{y}_1},\ldots,F_{\textbf{y}_s}$ and let the copula $C_H$ have continuous partial derivatives. Assume that $J$ is of bounded variation, continuous from above and with discontinuities of the first kind. Then
\begin{eqnarray*}
\lefteqn{\frac{1}{\sqrt{T}}\sum_{i=1}^T\left\{J(\widehat{F}_{\textbf{x}_1}(Z_{\textbf{x}_1}^{(i)}),\ldots,\widehat{F}_{\textbf{x}_k}(Z_{\textbf{x}_k}^{(i)}),\widehat{F}_{\textbf{y}_1}(Z_{\textbf{y}_1}^{(i)}),\ldots,\widehat{F}_{\textbf{y}_s}(Z_{\textbf{y}_s}^{(i)}))\right.}\\
&&-\left.E(J(F_{\textbf{x}_1}(Z_{\textbf{x}_1}^{(i)}),\ldots,F_{\textbf{x}_k}(Z_{\textbf{x}_k}^{(i)}),F_{\textbf{y}_1}(Z_{\textbf{y}_1}^{(i)}),\ldots,F_{\textbf{y}_s}(Z_{\textbf{y}_s}^{(i)})\right\}\\[0.3cm]
&\rightarrow &\int_{\left[0,1\right]^{k+s}}G(u_1,\ldots,u_{k+s})dJ(u_1,\ldots,u_{k+s}),
\end{eqnarray*}
where the limiting process and $G$ are Gaussian with zero mean.
\end{theorem}

In the following section we will conduct a simulation study of an M4 random field to assess the performance of the estimator
$\widehat{\widehat{\nu}}_{\mathbf{x},\mathbf{y}}^{\alpha,\beta}.$

\section{An M4 random field}

\pg It is well known that the class of max-stable processes called multivariate maxima of mo\-ving maxima processes or simply M4 processes, introduced by Smith and Weissman \cite{smi2}, is particularly well adapted to modeling the extreme bahaviour of several time series.

To illustrate the computation of the generalized madogram given in (\ref{mad_gen}) we will now define an M4 random field.

Lets consider that  the distribution of $(Z_{{\mathbi{x}}_1},\ldots,Z_{{\mathbi{x}}_p})$ is characterized by the copula
\begin{eqnarray}
C(u_{{\mathbi{x}}_1},\ldots,u_{{\mathbi{x}}_p})=\prod_{l=1}^{+\infty}\prod_{m=-\infty}^{+\infty}\bigwedge_{{\mathbi{x}}\in \{{\mathbi{x}}_1,\ldots,{\mathbi{x}}_p\}}u_{{\mathbi{x}}}^{a_{lm{\mathbi{x}}}},\quad u_{{\mathbi{x}}_i}\in [0,1],\ i=1,\ldots,p,\label{copula}
\end{eqnarray}
where, for each ${\mathbi{x}} \in \mathbb{Z}^2,$ $\{a_{lm{\mathbi{x}}}\}_{l\geq 1, m\in \mathbb{Z}}$ are non-negative constants such that $\displaystyle{\sum_{l=1}^{+\infty}\sum_{m=-\infty}^{+\infty}a_{lm{\mathbi{x}}}=1}.$
This random field ${\mathbf{Z}}$ is max-stable, since, for each $t>0,$ the copula (\ref{copula}) satisfies $$C^t(u_{{\mathbi{x}}_1},\ldots,u_{{\mathbi{x}}_p})=C(u^t_{{\mathbi{x}}_1},\ldots,u^t_{{\mathbi{x}}_p}),$$ for any locations ${\mathbi{x}}_1,\ldots,{\mathbi{x}}_p.$

As the M4 process considered in Smith and Weissman \cite{smi2}, we can consider that for each location ${\mathbi{x}},$  $Z_{{\mathbi{x}}}$ is a moving maxima of variables $X_{l,n},$ i.e.,
\begin{equation}
Z_{{\mathbi{x}}}=\max_{l\geq 1}\max_{-\infty<m<+\infty} a_{lm{\mathbi{x}}}X_{l,1-m},\quad {\mathbi{x}}\in \mathbb{Z}^2,\label{campo_M4}
\end{equation}
where $\{X_{l,n}\}_{l\geq 1, n\in \mathbb{Z}}$ is a family of independent unit Fréchet random variables. The dependence structure of $(Z_{{\mathbi{x}}_1},\ldots,Z_{{\mathbi{x}}_p})$ is regulated by the signatures patterns $a_{lm{\mathbi{x}}}$ and is given by (\ref{copula}).

For each pair of regions ${\mathbf{x}}=\{{\mathbi{x}}_1,\ldots,{\mathbi{x}}_k\}$ and ${\mathbf{y}}=\{{\mathbi{y}}_{k+1},\ldots,{\mathbi{y}}_{k+s}\}$ we have
\begin{eqnarray*}
V_{{\mathbf{x}},{\mathbf{y}}}(z_1,\ldots,z_k,z_{k+1},\ldots,z_{k+s})&=&-\ln C(e^{-z_1^{-1}},\ldots,e^{-z_{k+s}^{-1}})\\[0.3cm]
&=&
\sum_{l=1}^{+\infty}\sum_{m=-\infty}^{+\infty} \bigvee_{i=1}^{k+s}z_i^{-1}a_{lm\mathbi{x}_i},\ z_i\in \R,\ i=1,\ldots,k+s,
\end{eqnarray*}
and consequently, for $\alpha>0 $ and $\beta> 0$ we obtain
\begin{eqnarray*}
V_{{\mathbf{x}},{\mathbf{y}}}(\alpha,\ldots,\alpha,\beta,\ldots,\beta)=\sum_{l=1}^{+\infty}\sum_{m=-\infty}^{+\infty}\left( \bigvee_{i=1}^{k}\alpha^{-1}a_{lm\mathbi{x}_i} \vee \bigvee_{i=k+1}^{k+s}\beta ^{-1}a_{lm\mathbi{x}_i}\right).
\end{eqnarray*}

To illustrate the computation of the generalized madogram we shall consider, in what follows, examples with a finite number of signature patterns $(1\leq l\leq L)$ and a finite range of sequencial dependencies $(M_1\leq m\leq M_2)$.

\begin{ex}

{\rm{Lets consider that for each location ${\mathbi{x}}\in \mathbb{Z}^2$ with even coordinates we have \linebreak $a_{11{\mathbi{x}}}=a_{12{\mathbi{x}}}=\frac{1}{2}$ and otherwise $a_{11{\mathbi{x}}}=\frac{1}{4}=1-a_{12{\mathbi{x}}}.$ The values of
$(a_{11{\mathbi{x}}},a_{12{\mathbi{x}}})$ determine the moving pattern or signature pattern  of the random field, which in this case corresponds to one pattern ($L=1$).

\begin{figure}[!htb]
\begin{center}
\includegraphics[scale=0.45]{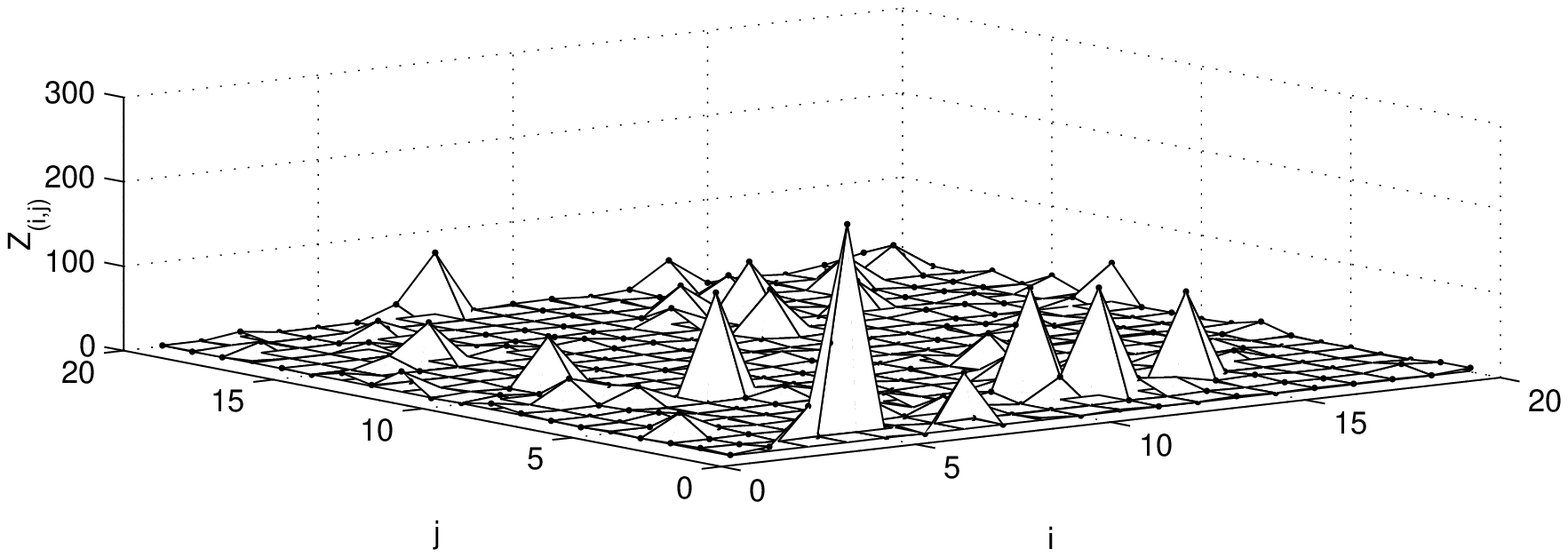}
\includegraphics[scale=0.45]{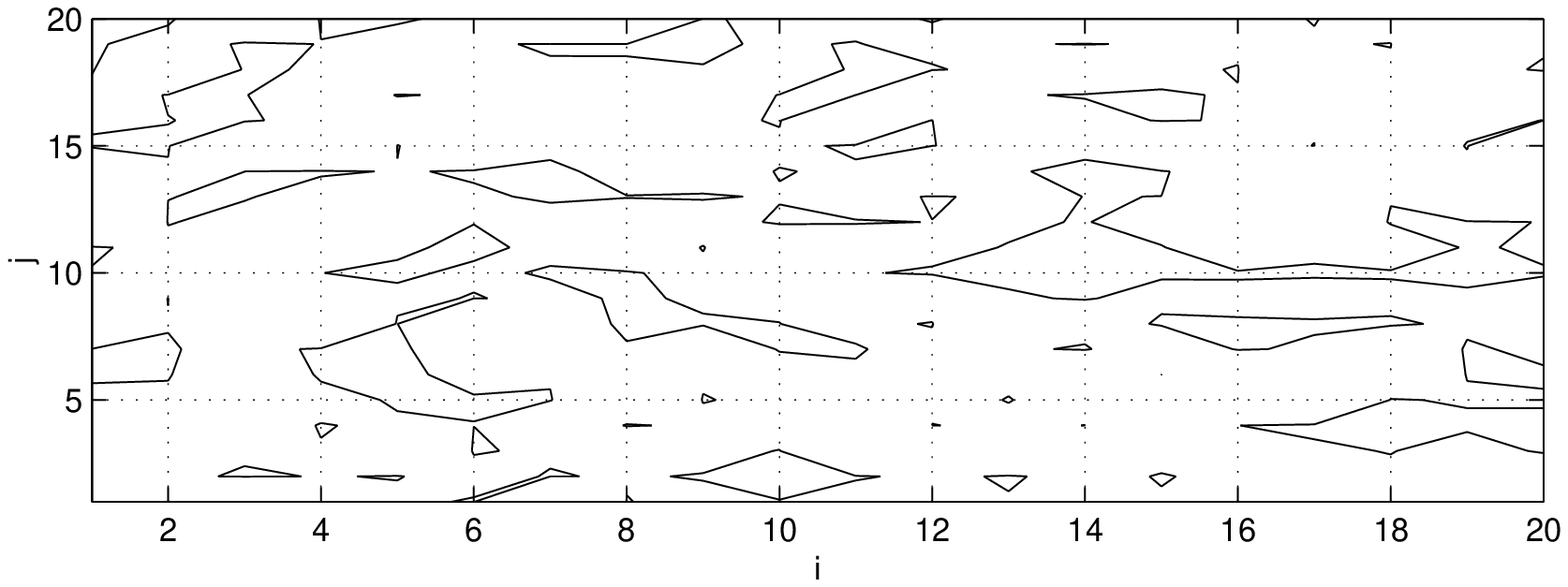}
\vspace{-0.5cm}\caption{\small Simulation of the M4 random field with $L=1$, $1\leq m\leq 2$ and for constants defined in Example 4.1.}
\end{center}
\label{fig:figura1}
\end{figure}

For the disjoint regions of locations ${\mathbf{x}}=\{(2,1),(2,2)\}$ and ${\mathbf{y}}=\{(3,3),(3,4)\}$ we have
$$V_{{\mathbf{x}},{\mathbf{y}}}(\alpha,\alpha,\beta,\beta)
= \frac{1}{4}(2\alpha^{-1} \vee \beta^{-1})+\frac{3}{4}(\alpha^{-1} \vee \beta^{-1})
$$ and therefore, the generalized madogram in this pair of locations is given by $$\nu^{\alpha,\beta}({\mathbf{x}},{\mathbf{y}})=\frac{\frac{1}{4}(2\alpha^{-1} \vee \beta^{-1})+\frac{3}{4}(\alpha^{-1} \vee \beta^{-1})}{1+\frac{1}{4}(2\alpha^{-1} \vee \beta^{-1})+\frac{3}{4}(\alpha^{-1} \vee \beta^{-1})}-\frac{1}{2}\left(\frac{\frac{5}{4}}{\alpha+\frac{5}{4}}+\frac{1}{\beta+1}\right),\quad \alpha>0,\ \beta>0.$$}}
\end{ex}

\begin{ex}
{\rm{Lets now assume that for each location ${\mathbi{x}}=(i,j)\in \mathbb{Z}^2,$ $a_{11{\mathbi{x}}}=\frac{1}{4}=1-a_{12{\mathbi{x}}}$ if $i\leq j$ and  $a_{11{\mathbi{x}}}=\frac{3}{4}=1-a_{12{\mathbi{x}}}$ if $i>j.$

As in the previous example the M4 random field generated by these sequences has a single signature pattern. Considering now two disjoint regions of locations with different size, ${\mathbf{x}}=\{(1,1)\}$ and ${\mathbf{y}}=\{(3,2),(3,3),(4,3)\},$  we obtain
$$V_{{\mathbf{x}},{\mathbf{y}}}(\alpha,\beta,\beta,\beta)
= \frac{1}{4}(\alpha^{-1}\vee 3\beta^{-1})+\frac{3}{4}(\alpha^{-1}\vee \beta^{-1})
$$
and consequently $$\nu^{\alpha,\beta}({\mathbf{x}},{\mathbf{y}})=\frac{\frac{1}{4}(\alpha^{-1}\vee 3\beta^{-1})+\frac{3}{4}(\alpha^{-1}\vee \beta^{-1})}{1+\frac{1}{4}(\alpha^{-1}\vee 3\beta^{-1})+\frac{3}{4}(\alpha^{-1}\vee \beta^{-1})}-\frac{1}{2}\left(\frac{1}{\alpha+1}+\frac{\frac{3}{2}}{\beta+\frac{3}{2}}\right),\quad \alpha>0,\ \beta>0.$$}}
\begin{figure}[!htb]
\begin{center}
\includegraphics[scale=0.45]{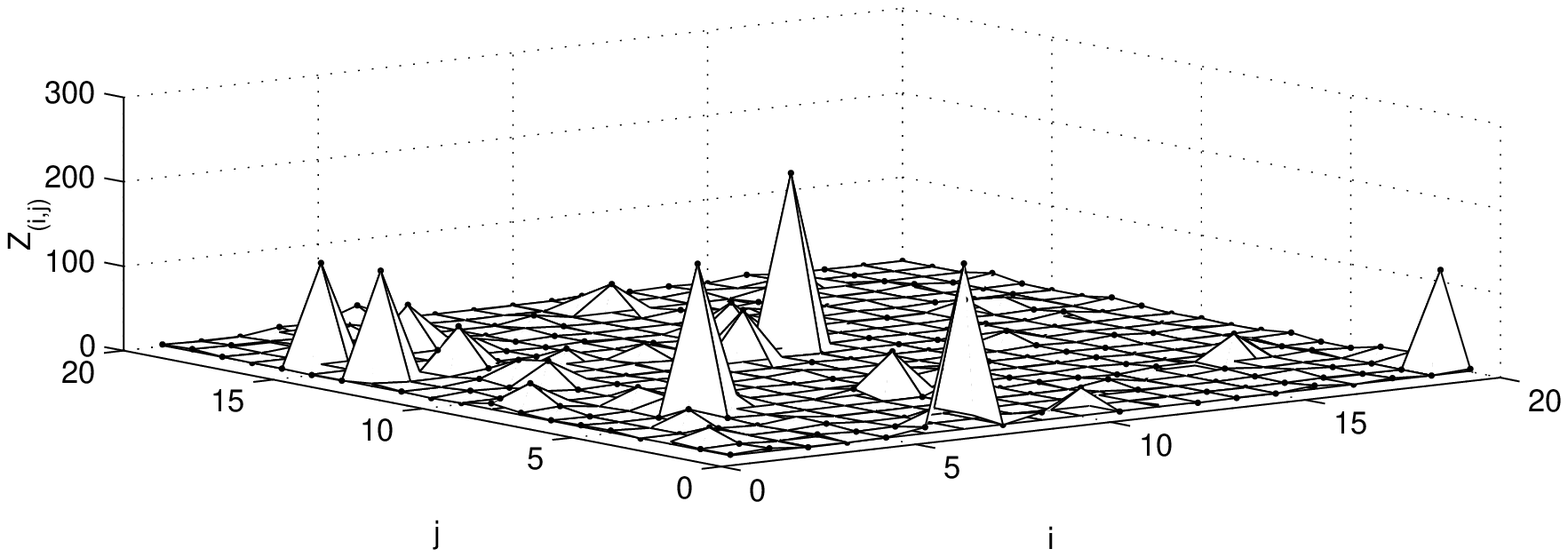}
\includegraphics[scale=0.45]{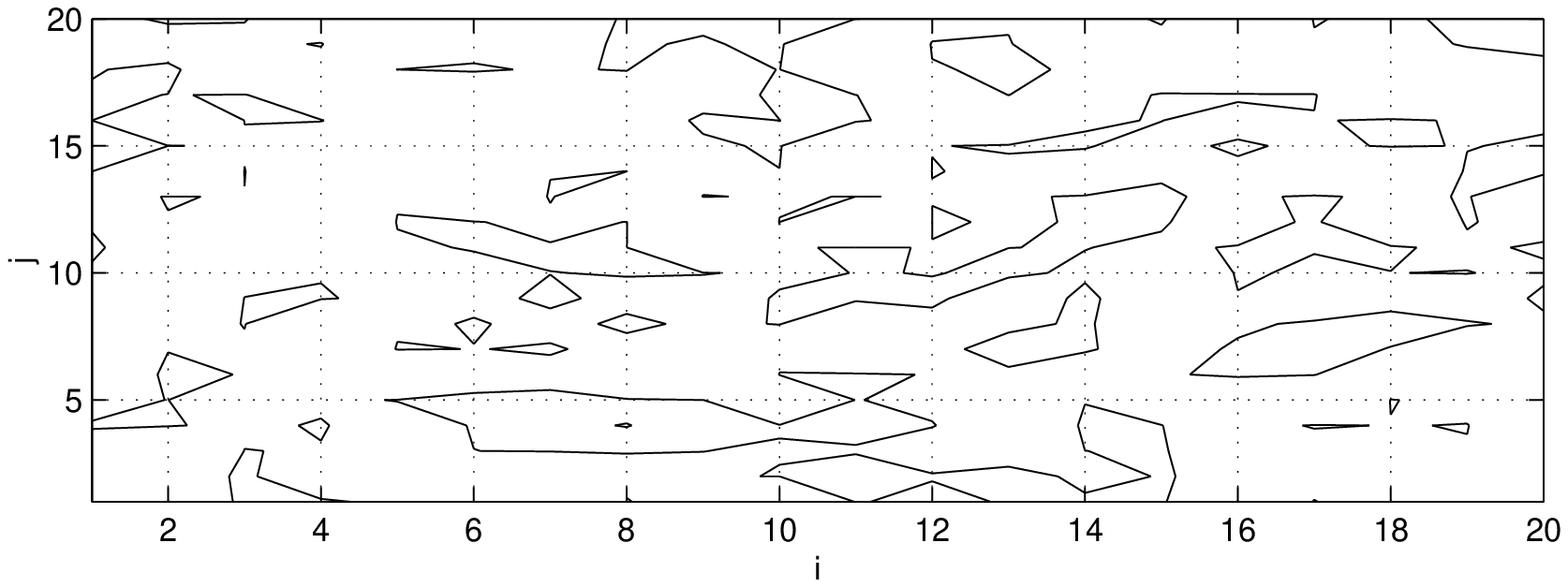}
\vspace{-0.5cm}\caption{\small Simulation of the M4 random field with $L=1$, $1\leq m\leq 2$ and for constants defined in Example 4.2.}
\end{center}
\label{fig:figura2}
\end{figure}

\end{ex}

\begin{ex}
{\rm{As stated in Zhang and Smith \cite{zhan1}, in a real data generating process it is unrealistic to assume that a single signature pattern would be sufficient to describe the shape of the process every time it exceeds some high threshold. Hence, we shall now consider one example with two signature patterns ($L=2$).

Lets assume that for each location ${\mathbi{x}}=(i,j)$ we have $a_{11{\mathbi{x}}}=a_{12{\mathbi{x}}}=a_{13{\mathbi{x}}}=\frac{1}{12},$ \linebreak $a_{21{\mathbi{x}}}=a_{22{\mathbi{x}}}=a_{23{\mathbi{x}}}=\frac{1}{4}$ if both coordinates are odd and $a_{11{\mathbi{x}}}=\frac{1}{18},\; a_{12{\mathbi{x}}}=\frac{1}{9},\; a_{13{\mathbi{x}}}=\frac{1}{6},$ \linebreak $a_{21{\mathbi{x}}}=a_{22{\mathbi{x}}}=a_{23{\mathbi{x}}}=\frac{2}{9}$ otherwise. Now the values of $(a_{11{\mathbi{x}}},a_{12{\mathbi{x}}},a_{13{\mathbi{x}}})$ and $(a_{21{\mathbi{x}}},a_{22{\mathbi{x}}},a_{23{\mathbi{x}}})$ define the two signature patterns of the random field.

\begin{figure}[!htb]
\begin{center}
\includegraphics[scale=0.45]{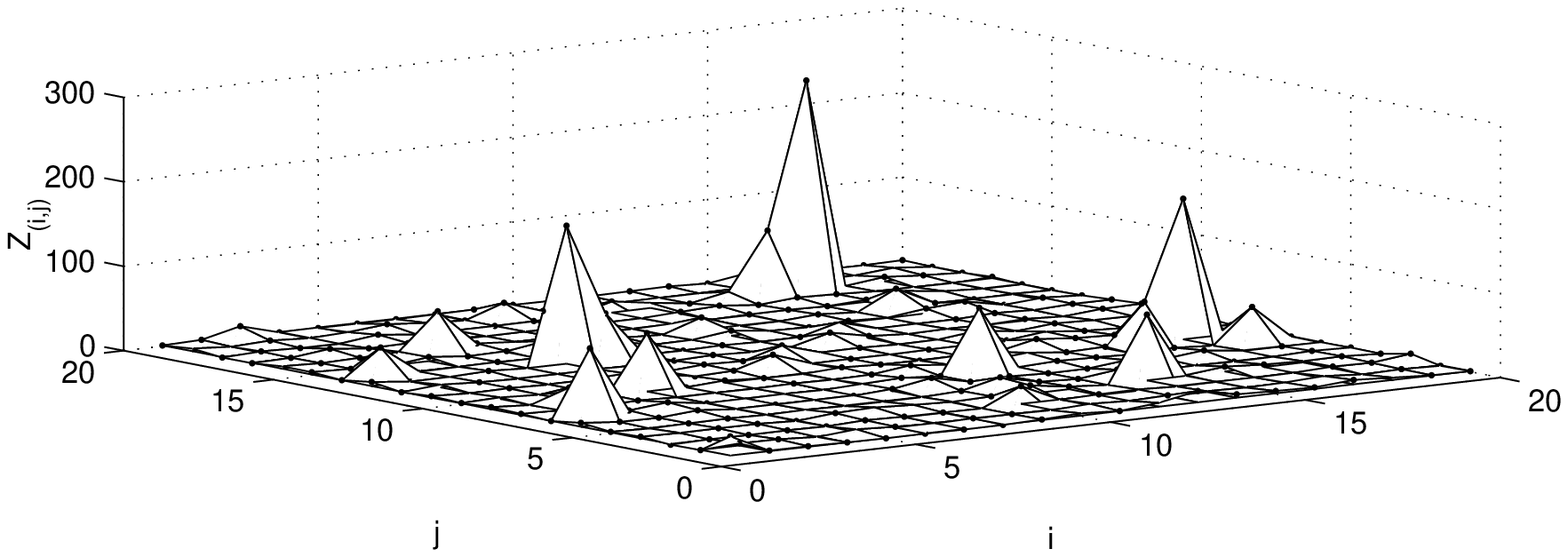}
\includegraphics[scale=0.45]{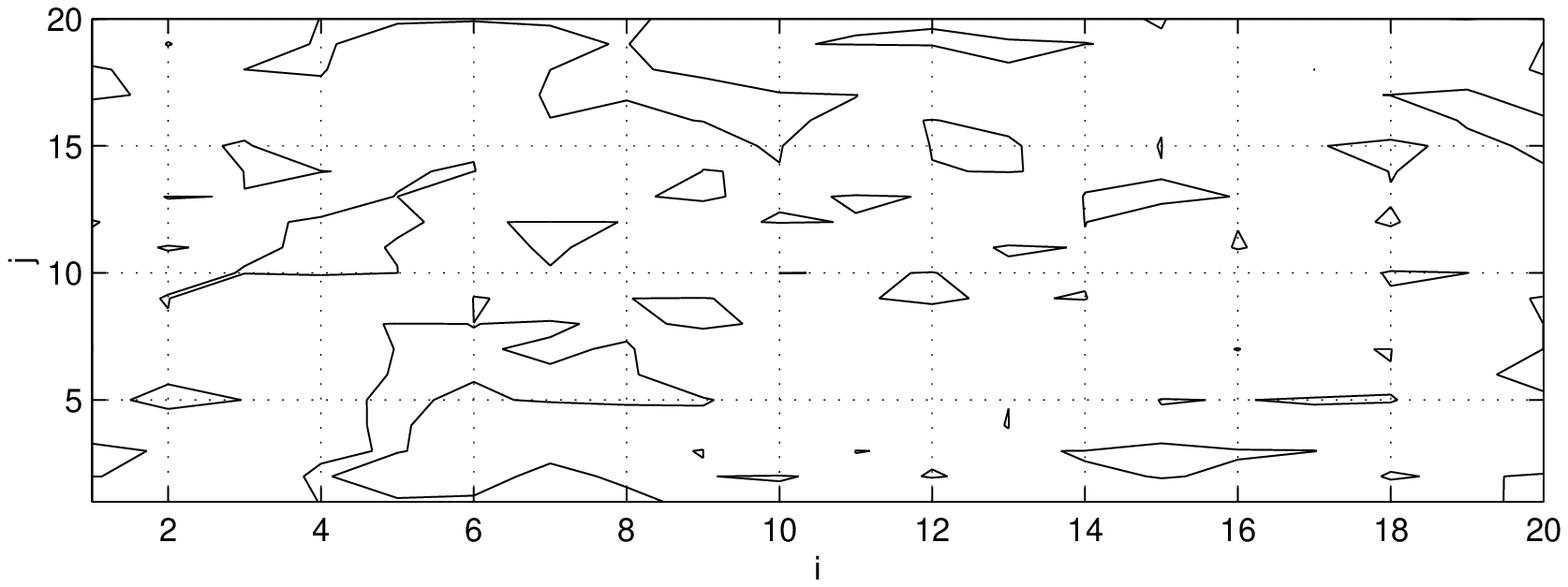}
\vspace{-0.5cm}\caption{\small Simulation of the M4 random field with $L=2$, $1\leq m\leq 3$ and for constants defined in Example 4.3.}
\end{center}
\label{fig:figura3}
\end{figure}

For two disjoint regions ${\mathbf{x}}=\{(2,1),(2,2)\}$ and ${\mathbf{y}}=\{(2,3),(3,3)\}$ we now have
$$V_{{\mathbf{x}},{\mathbf{y}}}(\alpha,\alpha,\beta,\beta)
= \left(\alpha^{-1}\frac{1}{18} \vee \beta^{-1}\frac{1}{12}\right)+\frac{1}{9}\left(\alpha^{-1} \vee \beta^{-1}\right)+\frac{1}{6}\left(\alpha^{-1} \vee \beta^{-1}\right)+\left(\alpha^{-1}\frac{2}{3} \vee \beta^{-1}\frac{3}{4}\right)
$$
and consequently
\begin{eqnarray*}
\nu^{\alpha,\beta}({\mathbf{x}},{\mathbf{y}})&=&\frac{\left(\frac{\alpha^{-1}}{18} \vee \frac{\beta^{-1}}{12}\right)+\frac{(\alpha^{-1} \vee \beta^{-1})}{9}+\frac{(\alpha^{-1} \vee \beta^{-1})}{6}+\left(\frac{2\alpha^{-1}}{3} \vee \frac{3\beta^{-1}}{4}\right)}{1+\left(\frac{\alpha^{-1}}{18} \vee \frac{\beta^{-1}}{12}\right)+\frac{(\alpha^{-1} \vee \beta^{-1})}{9}+\frac{(\alpha^{-1} \vee \beta^{-1})}{6}+\left(\frac{2\alpha^{-1}}{3} \vee \frac{3\beta^{-1}}{4}\right)} \\ && -\frac{1}{2}\left(\frac{1}{\alpha+1}+\frac{\frac{10}{9}}{\beta+\frac{10}{9}}\right),\quad \alpha>0,\ \beta>0.
\end{eqnarray*}
}}
\end{ex}

These examples will be used in the following simulation studies to assess the performance of the estimator given in (\ref{estim}).
Figures 4, 5 and 6  show the simulation results obtained by generating 50 replications of 100 independently and identically distributed max-stable M4 random fields in the three situations previously presented, with $\alpha$ and $\beta$ taking values in $\{k\times0.2:\;k=1,\ldots,100\}$.

As we can see from the values of the mean square error (Figures 4, 5 and 6) the estimates obtained from our estimator  $\widehat{\widehat{\nu}}^{\alpha,\beta}(\mathbf{x},\mathbf{y})$ are quite close to the true values of the generalized madogram.

\begin{figure}[!htb]
\begin{center}
{\tiny{\qquad\qquad$\nu^{\alpha,\beta}(\mathbf{x},\mathbf{y})\qquad\qquad\qquad \qquad\qquad \qquad\qquad\widehat{\widehat{\nu}}^{\alpha,\beta}(\mathbf{x},\mathbf{y})\qquad\qquad\;\qquad\qquad\qquad\qquad \qquad MSE\qquad\qquad$}}\\
\includegraphics[scale=0.23]{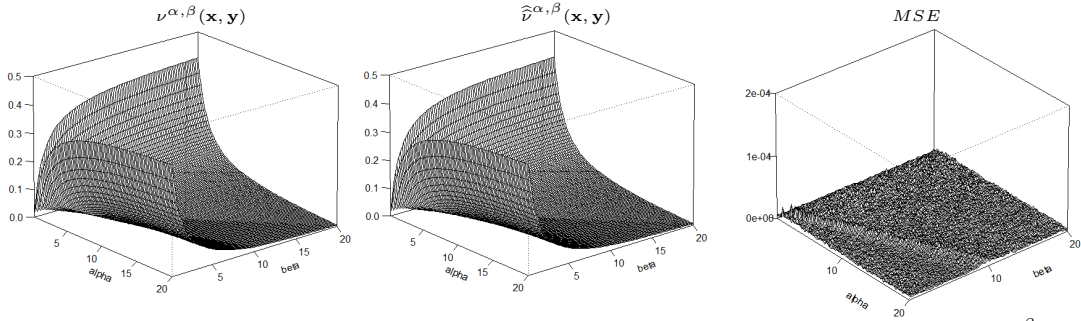}
\vspace{-0.5cm}\caption{\small The true values of the generalized madogram ($\nu^{\alpha,\beta}(\mathbf{x},\mathbf{y})$), the estimated values ($\widehat{\widehat{\nu}}^{\alpha,\beta}(\mathbf{x},\mathbf{y})$) and the estimated mean squared error ($MSE$) for  Example 4.1 (${\mathbf{x}}=\{(2,1),(2,2)\}$, ${\mathbf{y}}=\{(3,3),(3,4)\}$).}
\end{center}
\label{fig:figura4}
\end{figure}

\begin{figure}[!htb]
\begin{center}
{\tiny{\qquad\,$\nu^{\alpha,\beta}(\mathbf{x},\mathbf{y})\qquad\qquad\qquad\qquad\qquad\qquad \qquad\qquad\widehat{\widehat{\nu}}^{\alpha,\beta}(\mathbf{x},\mathbf{y})\qquad\qquad\qquad\qquad\qquad\qquad \qquad MSE$}}\\
\includegraphics[scale=0.23]{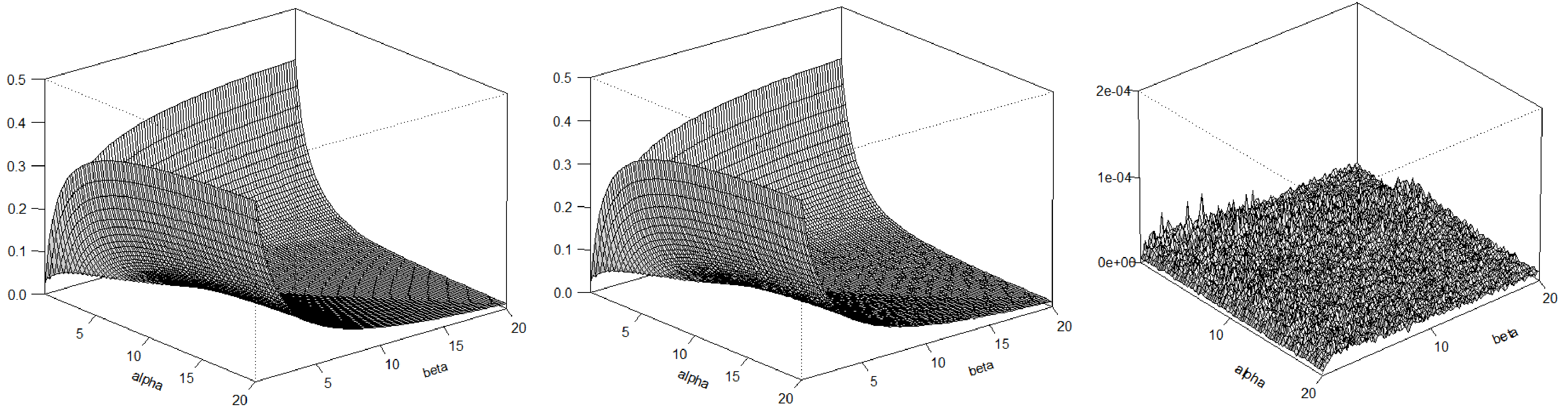}
\vspace{-0.5cm}\caption{\small The true values of the generalized madogram ($\nu^{\alpha,\beta}(\mathbf{x},\mathbf{y})$), the estimated values ($\widehat{\widehat{\nu}}^{\alpha,\beta}(\mathbf{x},\mathbf{y})$) and the estimated mean squared error ($MSE$) for  Example 4.2 (${\mathbf{x}}=\{(1,1)\}$, ${\mathbf{y}}=\{(3,2),(3,3),(4,3)\}$).}
\end{center}
\label{fig:figura5}
\end{figure}
\vspace{-0.7cm}
\begin{figure}[!htb]
\begin{center}
{\tiny{\qquad\,$\nu^{\alpha,\beta}(\mathbf{x},\mathbf{y})\qquad\qquad\qquad\qquad\qquad\qquad \qquad\qquad\widehat{\widehat{\nu}}^{\alpha,\beta}(\mathbf{x},\mathbf{y})\qquad\qquad\;\qquad\qquad\qquad\qquad \qquad MSE$}}\\
\includegraphics[scale=0.23]{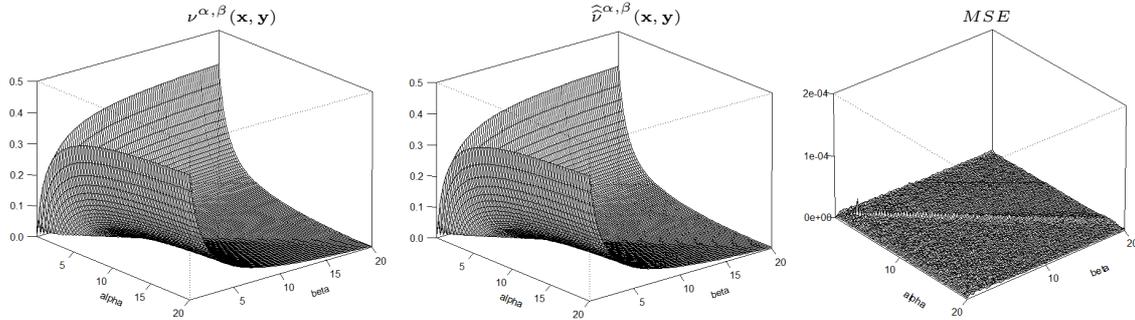}
\vspace{-0.5cm}\caption{\small The true values of the generalized madogram ($\nu^{\alpha,\beta}(\mathbf{x},\mathbf{y})$), the estimated values ($\widehat{\widehat{\nu}}^{\alpha,\beta}(\mathbf{x},\mathbf{y})$) and the estimated mean squared error ($MSE$) for Example 4.3 (${\mathbf{x}}=\{(2,1),(2,2)\}$, ${\mathbf{y}}=\{(2,3),(3,3)\}$).}
\end{center}
\label{fig:figura6}
\end{figure}

\section{Application to precipitation data}

\pg We now consider an application of the proposed estimator of the generalized  madogram to annual maxima values of daily maxima
rainfall in different topographic regions of Portugal.


In Portugal there are topographic differences between north and south, whereas in the north there are several mountains, in the south there are mainly plateaus and plains. The Central Cordillera, formed by the mountains of Sintra, Montejunto and Estrela, divides south and north and influences climate in Portugal, namely precipitation. The majority of the precipitation in Portugal comes from North-West and it is more abundant in the north than in the south due to the Central Cordillera that creates a physical barrier for precipitation.

\begin{figure}[!htb]
\begin{center}
\begin{minipage}[c]{0.3\linewidth}
\includegraphics[scale=0.5]{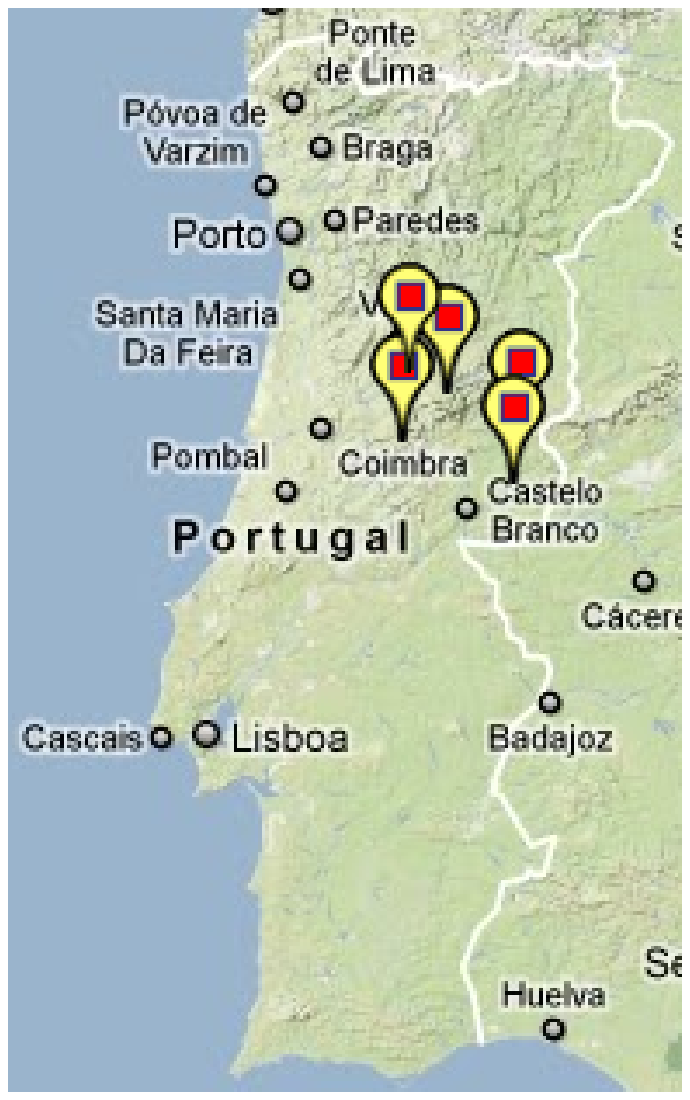}
\end{minipage}
\begin{minipage}[c]{0.3\linewidth}
\includegraphics[scale=0.4]{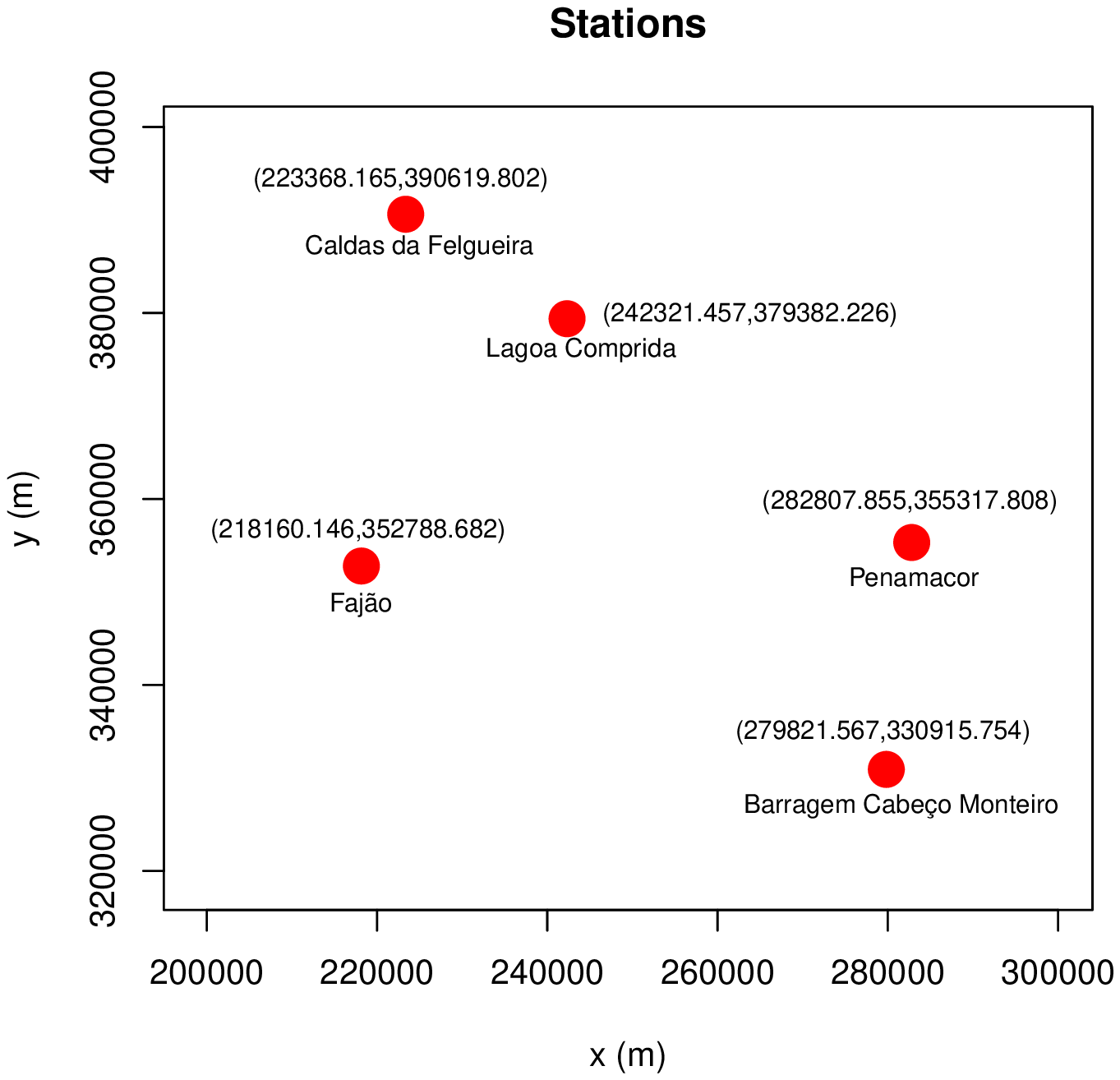}
\end{minipage}
\caption{The locations of the stations where precipitation data were collected, obtained from Portuguese National System of Water Resources (left) and their representation in Lambert coordinates (right).}
\label{fig:figura7}
\end{center}
\end{figure}

To study the dependence between extreme precipitation occurring in this mountain area and in surrounding areas, we considered precipitation data that consist of annual maxima of daily maxima  precipitation recorded over 38 years (1944-1981), in 5 Portuguese stations, obtained from the Portuguese National System of Water Resources (http://snirh.pt). In Figure \ref{fig:figura7} we can view the location of the 5 stations considered. We remark that the stations ``Lagoa Comprida'' and ``Fajão'' are located in North-West part of``Serra da Estrela'', the highest mountain in continental Portugal and part of the Central Cordillera.

 Since the data are maxima over a long period of time,  we assumed that they are independent over the years in each location. We also assumed that the random field is max-stable with unknown marginal distributions so data were previously transformed at each site so that they have a standard Fréchet distribution.

In Figure 8 we picture the estimated values of the generalized madogram, for several values of $\alpha$ and $\beta,$ when considering the mountain region $\mathbf{x}=\{\textrm{Fajão,\ L.Comprida}\}$ and different topographic surrounding regions, that are either to north or to the south of this region.


\begin{figure}[!htb]
\begin{center}
{\tiny{$\mathbf{x}=\{\textrm{L.Comprida,Fajão}\}\,\mathbf{y}=\{\textrm{B.C.M.}\}\,\,\,\,\,\,\mathbf{x}=\{\textrm{L.Comprida,Fajão}\}\,\mathbf{y}=\{\textrm {Penamacor}\}\,\,\,\,\,\,\mathbf{x}=\{\textrm{L.Comprida,Fajão}\}\,\mathbf{y}=\{\textrm{C.Felgueira}\}$}}
\includegraphics[scale=0.282]{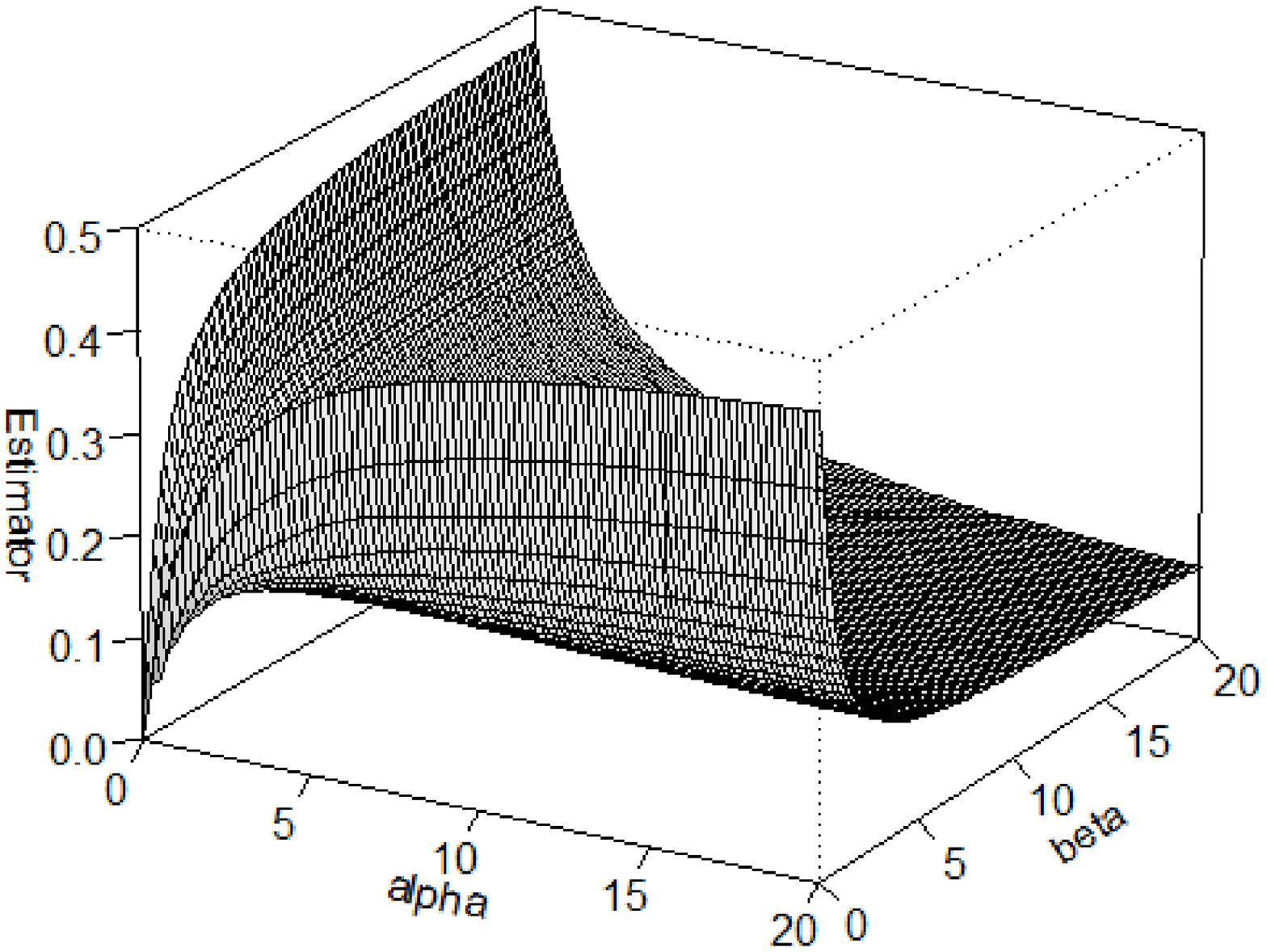}
\includegraphics[scale=0.282]{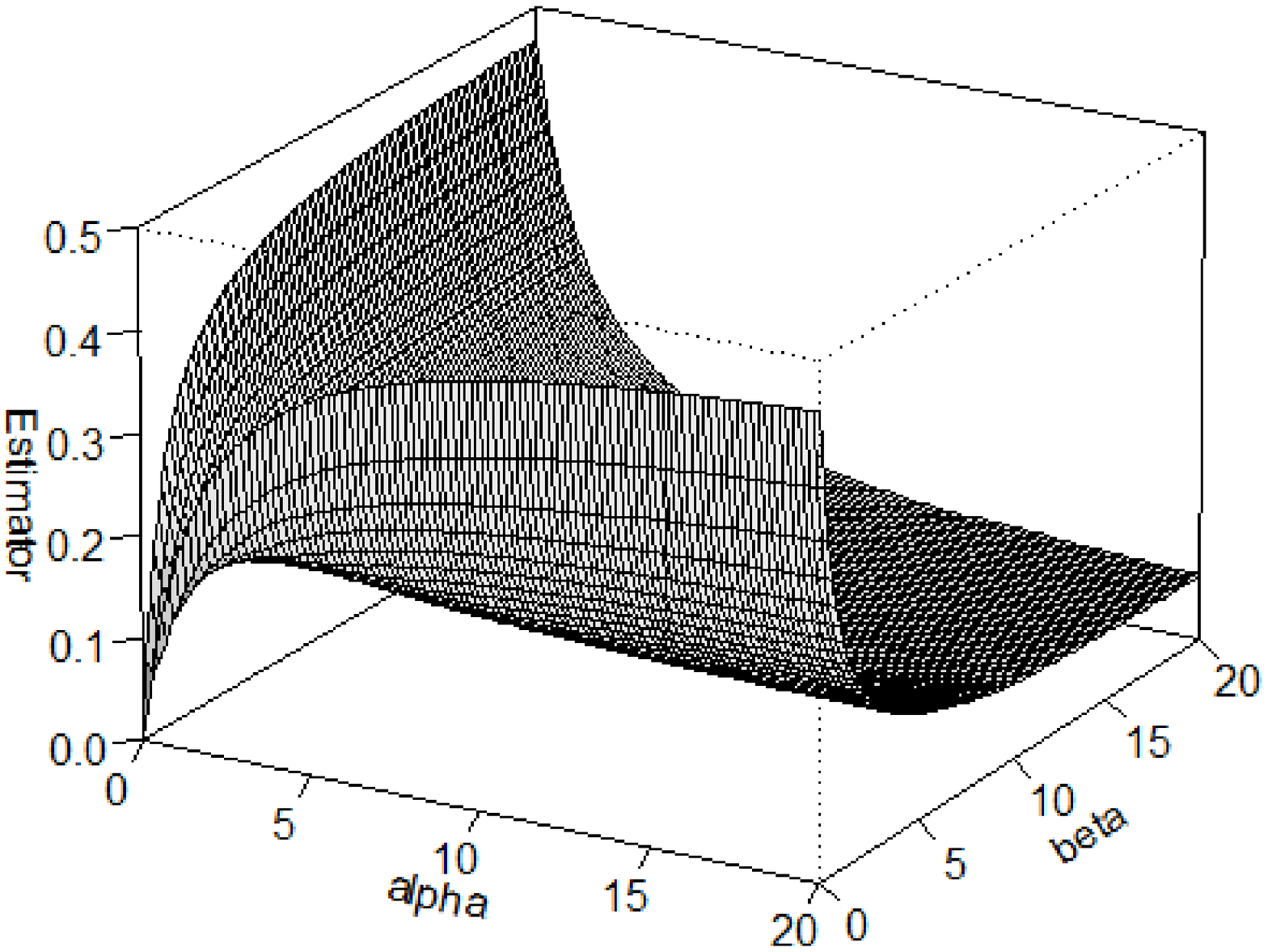}
\includegraphics[scale=0.282]{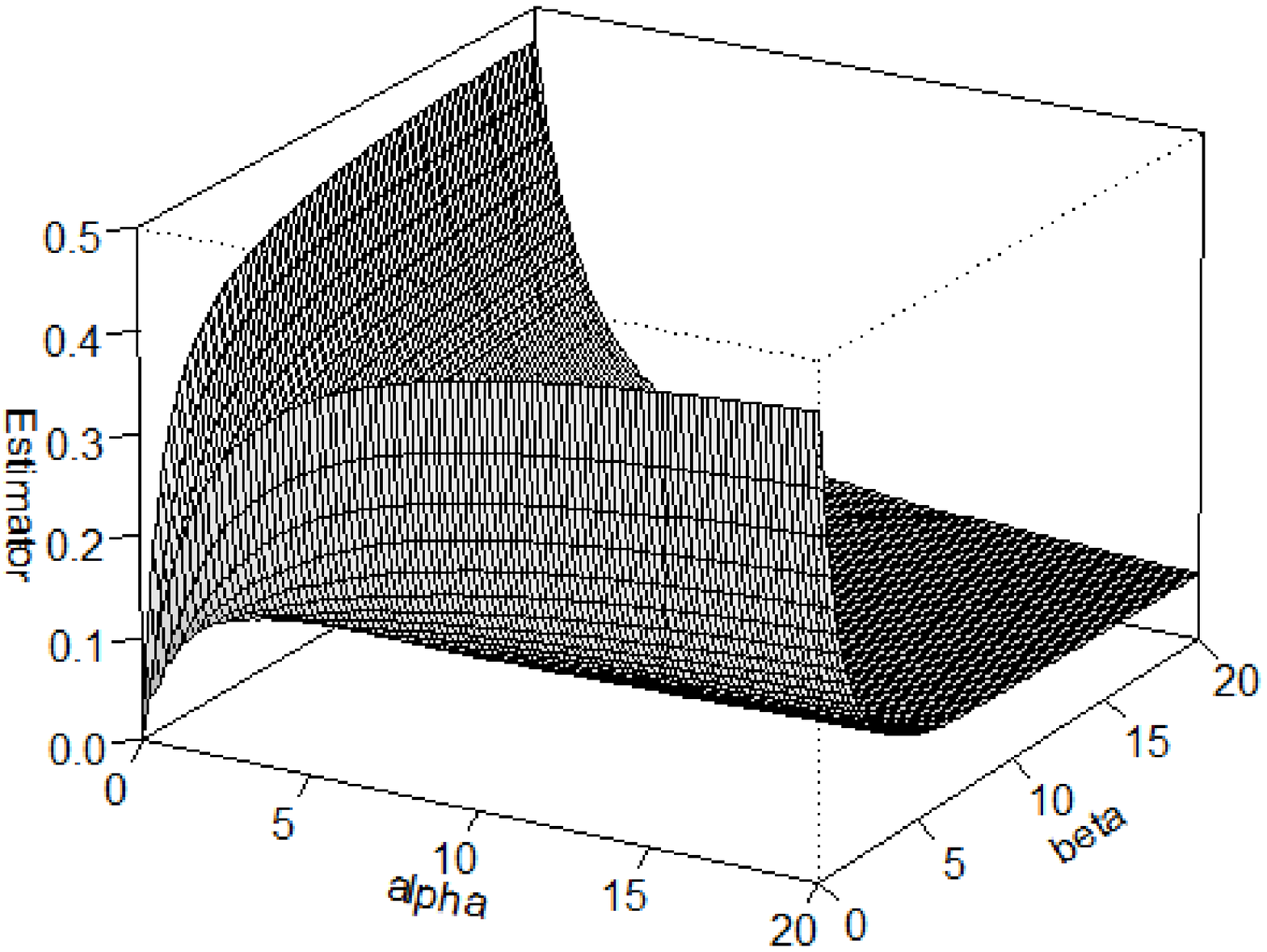}
\vspace{-0.2cm}\caption{\small  Generalized madogram ($\widehat{\widehat{\nu}}^{\alpha,\beta}(\mathbf{x},\mathbf{y})$) estimates obtained for each pair of disjoint regions and $\alpha$ and $\beta$  with values in $\{ k\times0.2:\;k=1,\ldots,100\}$.}
\end{center}
\label{fig:figura8}
\end{figure}

The estimated values for $\alpha=\lambda/k$ and $\beta=(1-\lambda)/s,$  $\lambda \in (0,1)$, where $k$ and $s$ are the number of locations in region $\mathbf{x}$ and $\mathbf{y}$, respectively, are presented in Figure 9, for the three situations considered in Figure 8.

\begin{figure}[!htb]
\begin{center}
{\tiny{$\qquad\qquad\qquad\widehat{\widehat{\nu}}^{\frac{\lambda}{2},1-\lambda}(\mathbf{x},\mathbf{y})\qquad\qquad$}}\\
\vspace{-0.8cm}
\includegraphics[scale=0.5]{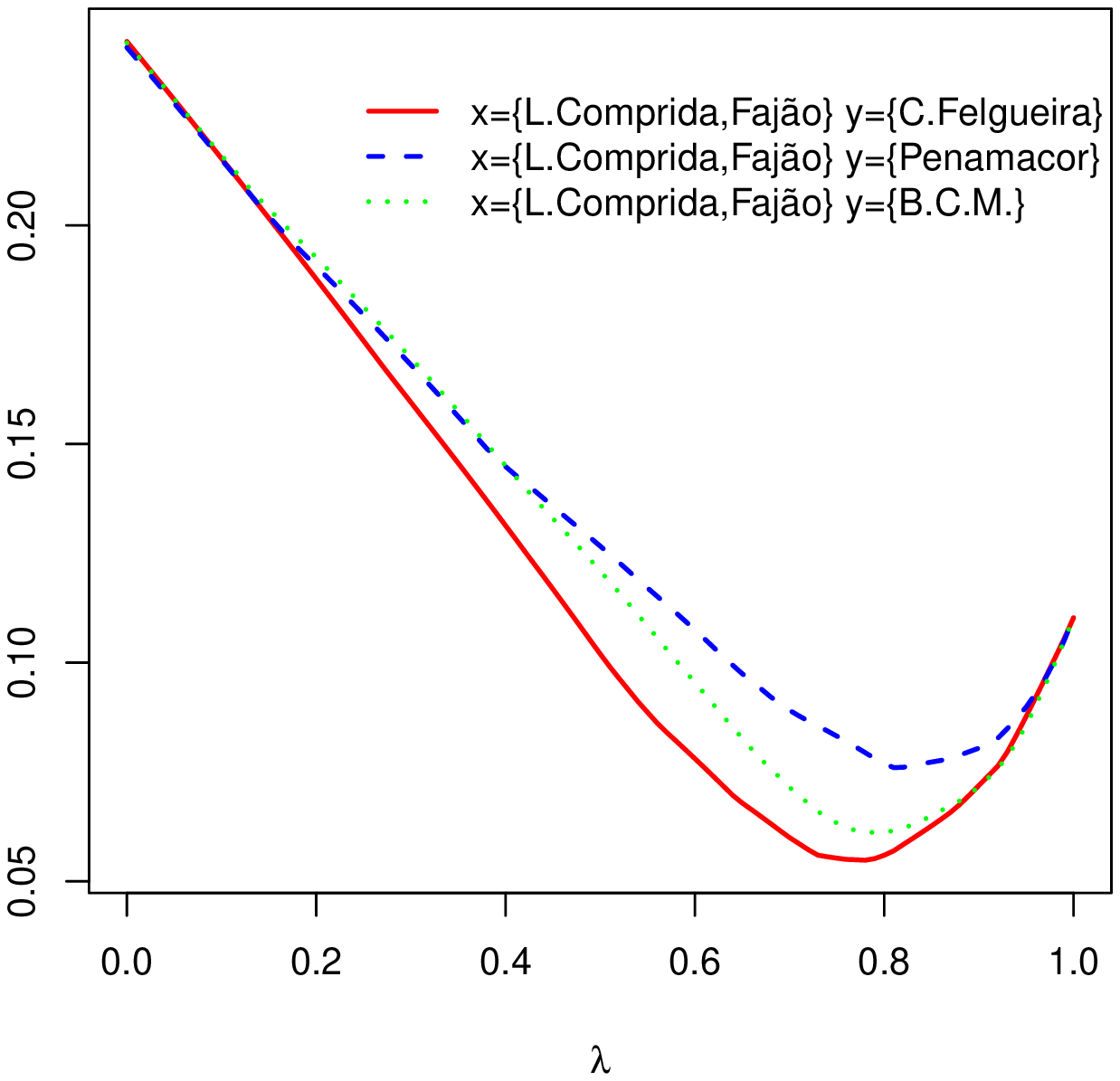}
\vspace{-0.5cm}\caption{\small  Generalized madogram ($\widehat{\widehat{\nu}}^{\alpha,\beta}(\mathbf{x},\mathbf{y})$) estimates obtained for each pair of disjoint regions when $\alpha=\frac{\lambda}{2}$ and $\beta=1-\lambda$ with $\lambda \in (0,1)$.}
\end{center}
\label{fig:figura9}
\end{figure}

Since lower values for $\nu^{\alpha,\beta}(\mathbf{x},\mathbf{y})$ indicate strong dependence, the results presented in Figure 9 suggest a stronger dependence between the mountain region $\mathbf{x}=\{\textrm{Fajão,\ L.Comprida}\}$ and the north region $\mathbf{y}=\{\textrm{C.Felgueira}\}$. This is in accordance with the previously stated that the Central Cordillera creates a physical barrier for precipitation in Portugal.

It would be interesting to further investigate this dependence with other regions but the lack of available data restricts the possible regions to study.

All the simulations presented in this paper were done in R statistical computing program \linebreak (http://cran.rproject.
org/). We remark that several packages on Extreme Value analysis have been recently introduced into R, but more recently
Ribatet \cite{rib} added to R the package \textsf{SpatialExtremes} that provides functions to analyze and fit max-stable processes to
spatial extremes.



\end{document}